\documentclass[12pt]{article}
\usepackage{amsmath,amsfonts,amssymb,amsthm}

\def\noi{\noindent}

\def\pf{\noi{\bf Proof.\ \,}}
\def\eop{{$\square$}}

\def\labtt#1{\label {#1} {\tt{<#1>}}}
\def\labttr#1{\label {#1} {\tt{<#1>}} \rm }

\def\U{\Upsilon}

\def\FF{{\mathbb F}}

\def\ZZ{{\mathbb Z}}

\def\la{\langle}
\def\ra{\rangle}
\def\<{\langle}
\def\>{\rangle}

\def\half{{1 \over 2}}

\def\bw#1{BW_{{2^{#1}}}}

\def\refpp#1{(\ref {#1})}
\begin{document}

\newtheorem{thm}{Theorem}[section]
\newtheorem{prop}[thm]{Proposition}
\newtheorem{lem}[thm]{Lemma}
\newtheorem{rem}[thm]{Remark}
\newtheorem{coro}[thm]{Corollary}
\newtheorem{conj}[thm]{Conjecture}
\newtheorem{de}[thm]{Definition}
\newtheorem{hyp}[thm]{Hypothesis}

\newtheorem{nota}[thm]{Notation}
\newtheorem{ex}[thm]{Example}
\newtheorem{proc}[thm]{Procedure}

\def\L{\Lambda}  
\begin{center}
{\Huge  Rank 72 high minimum norm lattices}

\vspace{10mm}
Robert L.~Griess Jr.
\\[0pt]
Department of Mathematics\\[0pt] University of Michigan\\[0pt]
Ann Arbor, MI 48109  \\[0pt]
USA  \\[0pt]
\vskip 1cm

\end{center} 

\newpage

\begin{abstract} 
Given a polarization of an even unimodular lattice and integer $k\ge 1$, 
we define a family of 
unimodular lattices $L(M,N,k)$.  Of special interest are certain $L(M,N,3)$ of 
rank 72.  Their minimum norms lie in $\{4, 6, 8\}$.  Norms 4 and 6 do occur.   
Consequently,  6 becomes the highest known minimum norm for rank 72 even unimodular lattices.  We discuss how norm 8 might occur for such a $L(M,N,3)$.  We note a few $L(M,N,k)$ in dimensions 96, 120 and 128 with moderately high minimum norms.  
\end{abstract}

\vskip 2cm


\newpage 

\tableofcontents

\section{Introduction}

Integral positive definite lattices with high norm for a given rank and discriminant have attracted a lot of attention, due to their connections with modular forms, number theory, combinatorics and group theory.  
Especially intriguing  are those even unimodular lattices which are {\it extremal}, i.e. their minimum norms  achieve the theoretical upper bound 
$2(\lfloor \frac n{24} \rfloor  +1)$, where $n$ is the rank.  
The rank  of an even unimodular lattices must be divisible by 8 (e.g., \cite{serre}).  The rank of an even integral unimodular extremal lattice is bounded (see \cite{bachocnebe} or Chapter 7 of \cite{cs} and the references therein).  
Extremal lattices are known to exist in dimensions a multiple of 8 up through 80, except for dimension 72.  An extremal rank 72 lattice would have   minimum norm  8  
\cite{cs,bachocnebe}.   

In this article, we construct a family of unimodular lattices $L(M,N,k)$ \refpp{gen1} for an integer $k$ and unimodular integral lattices $M, N$ which form a polarization \refpp{polarization2}.  Estimates on the minimum norm of $L(M,N,k)$ give 
some new examples of lattices with moderately high minimum norms.  

Of special interest are those $L(M,N,3)$ of dimension 72 where we input Niemeier lattices for $M$ and $N$.  Such a $L(M,N,3)$ have minimum norm 4, 6 or 8.  Norms 4 and 6 occur.  According to \cite{sloanesite}, our result is the first proof that 
there exists a rank 72 even unimodular lattice for which the minimum norm is at least 6.   
We indicate a specific criterion to be checked for such $L(M,N,3)$ to have minimum norm 8.    We conclude by noting certain $L(M,N,k)$ with moderately high norms  in dimensions 96, 120 and 128.

This work was supported in part by National Cheng Kung University where the author was a visiting distinguished professor; by 
Zhejiang University Center for Mathematical Research; by the 
University of Michigan; and by 
National Science Foundation Grant 
NSF (DMS-0600854).  We thank Alex Ryba for helpful discussions.  

\section{Integral sublattices of $\U^3$}

\begin{de}\labttr{lattice}  
In this article, {\it lattice} means a rational positive definite lattice.  The term {\it even lattice} means an integral lattice in which all norms are integral.  
For a lattice $L$, we define $\mu (L):=min \{ (x,x) \mid x\in L, x\ne 0 \}$ and call it the {\it minimum norm of $L$}.   If $L_1, L_2, \dots $ is a set of lattices, we define 
$\mu (L_1, L_2, \dots )$ to be the minimum of $\mu (L_1), \mu (L_2), \dots $.  
\end{de}

\begin{de}\labttr{polarization} Suppose that $E$ is an integral unimodular lattice.  A {\it polarization} is a pair of sublattices $X, Y$ such that 
$(X,X)\le 2\ZZ$, $(Y,Y)\le 2\ZZ$, $X+Y=E$ and $X\cap Y =2E$.  It follows that $E$ is even.  
If $E$ is a lattice and $r>0$ is a rational number such that 
$\sqrt r \, E$ is an integral unimodular lattice, a 
{\it polarization} of $E$ is a pair of sublattices $X, Y$ so that 
$\sqrt rM, \sqrt rN$ is a polarization of $\sqrt rE$.  
\end{de}  

\begin{rem}\labttr{polarization2} If $Z$ is one of $X, Y$ as in \refpp{polarization} 
and $E$ is unimodular, then $\frac 1{\sqrt 2} \, Z$ is integral and unimodular, but may not be even.  If $\frac 1{\sqrt {2}} X$ and $\frac 1{\sqrt {2}} Y$ are both even lattices we call the polarization an {\it even polarization}.   If $E$ is not unimodular but $\sqrt r E$ is, the polarization $X, Y$ of $E$ is called {\it even} if the polarization $\sqrt r X, \sqrt r Y$ is even.  
\end{rem} 

\begin{nota}\labttr{ups} 
We let $\U$ be a lattice so that $U:= \sqrt 2 \,  \U$ is an even, integral unimodular lattice.  
\end{nota}

A polarization of $\U$ is therefore a pair of integral sublattices $M, N$ such that 
$M+N=\U$ and $M\cap N = 2\U$.  

For the time being, $rank(\U )=rank(U)$ is an arbitrary multiple of 8.    
We know the complete list of possibilities for even, integral unimodular lattices only in dimensions 8, 16 and 24.  The rank 24 lattices are called {\it Niemeier lattices} since they were first classified by Niemeier \cite{niem}.  

\begin{lem}\labtt{e8polar}  The $E_8$-lattice has an even polarization.  
\end{lem}
\pf  This is a standard fact.  It follows since the $E_8$ lattice modulo 2 has a nonsingular form with maximal Witt index. One then quotes the characterization of $E_8$ as the unique (up to isometry) rank 8 even unimodular lattice.   Another proof uses the existence of a fourvolution \refpp{fourvolution} on $E_8$ (one exists, for example, in a natural $Weyl(D_8)$ subgroup; if one identifies $E_8$ with $\bw 3$, the natural group of isometries $\bw 3$ contains lower fourvolutions).  
\eop

\begin{nota}\labttr{gen1}
We use the notation of \refpp{ups} and let $M, N$ be a polarization of $\Upsilon$.  
Let $k\ge 2$.  Define these sublattices of $\U^k$: 
$$L_M:=\{ (x_1, \dots , x_k)\in M^k \mid x_1+\cdots + x_k \in M\cap N \},$$ 
$$L^N:= \{(y,y,\dots , y) \mid y \in N \},$$  
$$L(M,N,k):=L_M + L^N.$$  
\end{nota}

\begin{rem}\labttr{gen1.5} Because $L(M,N,1)=N$ and $L(M,N,2)\cong U\perp U$, the interesting case is $k \ge 3$.  If $k=2q$ is even, $L(M,N,k)$ contains 
$L^M +L^N$, 
a sublattice isometric to 
$\sqrt {q}\, U$ 
\end{rem}

\begin{prop}\labtt{gen2} 
(i) The lattice $L(M,N,k)$ is an integral lattice and the sublattice $L_M$ is even.  

(ii) If $k$ is an even integer or $N$ is an even lattice, 
$L(M,N,k)$ is an even lattice.  Otherwise, $L(M,N,k)$ is odd.  

(iii) $L(M,N,k)$ is unimodular.  
\end{prop}
\pf   
(i) 
To prove integrality, 
one shows that $L_M$ and $L^N$ are integral lattices and that 
$(L_M,L^N)\le \ZZ$.  The latter follows since 
for $(x_1, \dots , x_k) \in L_M$, $\sum_i x_i \in N$, an integral lattice. 
Finally, the evenness of $L_M$ is obvious since it is integral and a set of generators is even (e.g., all vectors of the form $(x,x,0^{k-1}), x\in M$ and $(y,0^{k-1}), y\in 2\U$).  

(ii) This is obvious from the  definition of $L^N$.  

(iii) 
To prove unimodularity, it suffices by \refpp{indexdet}  to show that $|L: L_M|^2 = det(L_M)$.   
We have $det(L_M)=det(M^k)|M^k:L_M|^2=1\cdot 2^{rank(M)}$ and 
$|L:L_M|=|L_M+L^N:L_M|=|L^N : L^N\cap L_M|=|L^N:L^N\cap M^k||L^N\cap M^k : L^N\cap L_M|=2^{\half rank(M)}\cdot 1$.   
\eop

\begin{thm}\labtt{minlmn}  We use the notation $\mu (L_1, L_2, \dots )$ \refpp{lattice}.  

(i) $\mu (L_M)=2\mu (M,U)$ and $\mu (L^N)=k\mu (N)$.  

(ii) $\mu (L) \le min\{ k\mu (N), 2\mu (M,U) \}$.  

(iii) 
$\mu (L) \ge min\{ \frac k2 \mu (U), 2\mu (M,U)\}$.    
\end{thm} 
\pf  
(i) To determine $\mu (L_M)$, consider the possibility that all entries of $(x_1, \dots, x_k) \in L_M$ are in $2\U$.

(ii) This follows from (i) since $L_M$ and $L^N$ are sublattices of $L$.

(iii) If a vector is in $L\setminus L_M$, all of its coordinates are nonzero.  
\eop

\begin{nota}\labttr{leechdef} We let $\L$ be a {\it Leech lattice}, i.e., a Niemeier lattice without roots.  
\end{nota}

Uniqueness of a rootless Niemeier lattice  was proved first in \cite{conway}, then in different styles in \cite{borcherds} and  \cite{poe}.

 We illustrate the use of \refpp{minlmn} by constructing a Leech lattice.  This argument comes from 
 \cite{tits}, \cite{lm}.  
An analogous construction of a Golay code was created earlier by Turyn \cite{turyn}.  
The original existence proof of the Leech lattice \cite{leech} makes use of the Golay code (whereas \refpp{leech} does not). 

\begin{coro}\labtt{leech}
Leech lattices exist.  
\end{coro}
\pf   
We take $M\cong N \cong E_8$ \refpp{e8polar}.  
From \refpp{minlmn},  $3 \le \mu (L)\le 4$.  Since $L(M,N,3)$ is even, $\mu (L(M,N,3)=4$.
\eop 

\begin{nota}\labttr{leechnota} 
We use the standard notation $\L$ for a Leech lattice.  
\end{nota}

\section{Minimum norms for rank 72 $L(M,N,3)$}

\begin{nota}\labttr{rank72nota} 
In this section, $L(M,N,3)$ is a rank 72 lattice for which $M$ and $N$ are Niemeier lattices.  
\end{nota} 

The minimum norm of a Niemeier lattice is 2 unless it is the Leech lattice, for which the minimum norm is 4.   

\begin{coro}\labtt{mul72}
(i) 
$\mu (L(M,N,3)) \ge 4$.  

(ii) If $M \not \cong \L$, then $\mu (L(M,N,3)) = 4$.  

(iii) If $U\cong M\cong \L$, then 
$\mu (L(M,N,3)) \ge 6$.  

(iv) If $U\cong M\cong \L$, 
and $N\not \cong \L$, 
then 
$\mu (L(M,N,3)) = 6$.  
\end{coro} 

We now prove that 
situations (ii) and (iv) of the Corollary actually occur.  This means proof that suitable polarizations of $\U$ exist.  

\begin{prop}\labtt{n4n6}
There exist $L(M,N,3)$ 
with minimum norms 4 and 6. 
\end{prop}
\pf 
We take $U\cong E_8^3$ and $M, N \le U, M\cong N\cong \sqrt 2 E_8^3$ such that 
$M+N=U$ (for example, the orthogonal direct sum of three polarizations as in \refpp{leech} will do).  
Then  (ii) applies. 

If $U\cong \L$, take in $\U$ any sublattice $M\cong \L$ (see \refpp{fourvolutiontype}, \refpp{leechleechpolar}) and any $N\cong E_8^3$ (see \cite {poe} for existence).  Then (iv) applies.  
\eop

\begin{coro}\labtt{n8?} 
If $\mu (L(M,N,3))=8$,  $M\cong N\cong \L$.  
\end{coro}

 The question remains whether there exists a polarization $M, N$ so that $\mu (L(M,N,3))=8$.

\begin{rem}\labttr{niemniem} 
It would be useful to know more about embeddings of 
$\sqrt 2 J$ into $K$, where $J, K$ are Niemeier lattices.  For the case $K\cong \L$, see \cite{dmn}, Th. 4.1.  Note also that embeddings of 
$\sqrt 2 E_8^3$ in $\L$ were used extensively in \cite{poe}.  
\end{rem}

\section{Norm 6 vectors in rank 72 $L(M,N,3)$ }

\begin{nota}\labttr{norm6nota} 
Let $L:=L(M,N,3)$, where 
$M\cong N\cong \L$ (by \refpp{leechleechpolar}, there exists such a polarization).  
\end{nota} 

From \refpp{mul72}(iii), $\mu(L)\ge 6$.  
We consider the possibility that $L$ has vectors of norm 6 and derive some results about forms of norm 6 vectors.   

We use parentheses both for inner products $(x,y)$ and $n$-tuples $(x_1, \dots, x_n)$.  We hope for no confusion when $n=2$.   

\begin{nota}\labttr{setup2} 
We call an ordered 4-tuple $(w,x,y,z)\in N\times M\times M\times M$ {\it admissible} if $x+y+z\in M\cap N$.  
The elements of $L$ are the $(x+w,y+w,z+w)$, for all 
 admissible 4-tuples $(w,x,y,z)$.  
 We call admissible 4-tuples 
 $(x,y,z,w)$  and $(x',y',z',w')$ 
 {\it equivalent} if 
$(x+w,y+w,z+w)=(x'+w',y'+w',z'+w')$.  
An {\it offender} is a 4-tuple $(x,y,z,w)$ 
such that  each of $r_x:=x+w, r_y:=y+w, r_z:=z+w$ has norm 2. 
Offenders are those admissible 4-tuples which give norm 6 vectors $(x+w,y+w,z+w) \in L$ (since $\mu (M)=4$, $w\notin M$ or else $M$ would contain roots).   The set $r_x, r_y, r_z$ is called {\it a triple of offender roots}.  
\end{nota}

If there are no offenders, $L$ has minimum norm 8.  We therefore study hypothetical offenders.    

The rational lattice $\U = M+N$ is not integral (in fact, $(\U, \U)=\half \ZZ$).   
The next result asserts integrality of the sublattice of $\U$  spanned by  the components of an offender. 

\begin{lem}\labtt{offint}
For an offender, $(w,x,y,z)$, we define $K$ to be the $\ZZ$-span of $w,x,y, z$.  
Then 

(i) The image of $K$ in $(M+N)/M$ has order 2;

(ii) $K$ is an even integral lattice. 
\end{lem}  
\pf
(i) The image of $K$ in $(M+N)/M$ is spanned by the image of $w$, and $w\notin M, 2w\in M$.  

(ii) Since $x, y, z$ lie in an integral lattice $M$ and $w\in N$ is integral, it suffices to prove that each of $(w,x), (w,y), (w,z)$ is integral.  We have $2=(w+x,w+x)=(w,w)+2(w,x)+(x,x)$.  Since  $M$ and $N$ are even lattices, $(w,w)$ and $(x,x)$ are even integers.  So $(w,x)$ is integral.  Similarly, we prove $(w,y), (w,z)$ are integral. 
\eop 

\begin{lem}\labtt{shortmod} 
Let $Q$ be a sublattice of $\L$, $Q\cong \sqrt 2 \L$.  
The $2^{12}-1$ nontrivial cosets each contain exactly 48 norm 4 vectors, and such a set of 48 is an orthogonal frame: two members are proportional or orthogonal.  
\end{lem}
\pf
This may be proved by a rescaling of the argument that in $\L$, the  norm 8 vectors which lie in the same coset of $2\L$ 
constitute an orthogonal frame of 48 vectors.  See \cite{conway,gr12}.  
\eop 

\begin{lem}\labtt{wnorm4} Suppose that $M$ has fourvolution type \refpp{fourvolutiontype}.  
If $(w,x,y,z)$ is admissible and $w\notin M$, there exists an equivalent 
admissible quadruple  $(w',x',y',z')$  such that $w'$ has norm 4. 
\end{lem}
\pf
This follows from \refpp{shortmod}.  There exists $v\in \U$ so that 
$w':=w-2v\in N$ has norm 4 (recall that $2\U=M\cap N$).  Take $x':=x+2v, y':=y+2v, z':=z+2v$.   These three vectors lie in $M$.   
\eop

\begin{lem}\labtt{orthogoffenderroots}  
A triple of offender roots is a pairwise orthogonal set.  
\end{lem}  
\pf
Suppose that two such roots are not orthogonal, say $r=w+x$ and $s=w+y$.  
Define $J:=span\{r, s\}$, an $A_2$-lattice (note that $J$ is integral, by \refpp{offint}(ii)).  
Since  $M\cap J$  is contained in $M$, it is rootless.  
However,  $M\cap J$ has index 2 in $J$ gives a contradiction since every index 2 sublattice of $J$ contains roots.  
\eop

\begin{lem} \labtt{ipseq}
Let $r, s, t$ be the three roots from an offender triple (in any order).  
The unordered set of inner products $(w,r), (w,s), (w,t)$ is $0,0,\pm1$.  
The unordered set of norms for $x, y, z$ is one of  $6, 6, 4$ or $6, 6, 8$.  
\end{lem}
\pf  
The second statement follows from the first, which we now prove.  
Let $r'\in \{r,-r\}$ satisfy $(w,r')\le 0$.  Similarly, let 
$s'\in \{s,-s\}$ satisfy $(w,s')\le 0$
and 
$t'\in \{t,-t\}$ satisfy $(w,t')\le 0$. 
Then $w+r'+s'+t'\in M\cap N$ and $w+r'+s'+t'$ has norm $4+2+2+2+e$, where 
$e\le 0$ and $e$ is even.  

We observe that if $w+r'+s'+t'$ were 0, the pairwise orthogonality of $r, s, t$ would imply that $w$ has norm 6, which is not the case.  Therefore, 
$w+r'+s'+t'$ has even norm at least 8.  Consequently, 
$e=0$ or $e=-2$.   Since $M\cap N\cong \sqrt 2 \L$, in which norms are divisible by 4 
and nonzero norms are at least 8, 
$e=-2$.  Therefore all but one of $(w,r), (w,s), (w,t)$ is 0 and the remaining one is $\pm 1$.  
\eop

\begin{nota}\labttr{super } 
An offender $(w,x,y,z)$ is a   {\it super  offender}  if $w$ has norm 4 and the norms of $x, y, z$ in some order are 6, 6, 4. 
\end{nota}

\begin{lem}\labtt{44}
We may assume that an offender $(w,x,y,z)$ satisfies $(w,w)=4$, $(w,t)=1$ and $(z,z)=4$.  In other words, if an offender exists,  a super offender exists.  
\end{lem}
\pf 
Since $(w,t)=\pm 1$, $z=t-w$ has norm 4 or 8, respectively. 
Suppose the latter.  Then  $(-w,-x,-y, z+2w)$ is admissible and its final component $z+2w=t+w$ has norm 4.   
Therefore, $(-w,-x,-y, z+2w)$ is a super offender.  
\eop

\begin{thm}\labtt{6or8}  Let $L:=L(M, N)$, where $M\cong N$ are isometric to the Leech lattice.  Then the minimum norm of $L$ is 6 if and only if there exists a super offender.  Otherwise, the minimum norm is 8.  
\end{thm}

\begin{rem}\labttr{conclusion}
Given $M, N$, \refpp{6or8} indicates that checking  
a (very large) finite number of inner products 
will settle $\mu (L(M,N,3))$.  

There are finitely many polarizations $M, N$ of $\U$.  Possibly some 
$L(M,N,3)$ 
have minimum norm 6 and others have minimum norm 8.  

Use of isometry groups and other theory might reduce the number of  computations significantly.  
\end{rem}

\newpage 

\section{Some higher dimensionss}

\begin{lem}\labtt{gen3.5}  
There exist rank 32 even integral unimodular lattices $U, M, N$ so that 
$\mu (U)=\mu(M)=4$, $\mu(N)\in \{2,4\}$ and 
$\sqrt 2 M, \sqrt 2 N$ is a polarization of $U$.  
\end{lem}  
\pf
 We take $U$ to be $\bw 5$. If $f$ is a fourvolution in $O(U)$, then 
$M:=(f-1)U\cong \sqrt 2 U$.  
Therefore, the natural $\FF_2$-valued quadratic form on $U/2U$ is split (i.e.,  has maximal Witt index) and so there exists an even unimodular lattice $N$ so that $\sqrt 2 N$ is between $U$ and $2U$ and $\sqrt 2 N/2U$ complements $M/2U$ in $U/2U$.  
The extremal bound $\mu (N)\le 4$ and evenness of $N$ imply the last statement.  
\eop 

\smallskip

We now exhibit a few even unimodular lattices for which the minimum norm is moderately  close to the extremal bound 
$2(1+\lfloor \frac {rank(L)}{24}
\rfloor)$.

\begin{prop}\labtt{gen4}  Let $U, M, N$ be as in \refpp{gen3.5} and let $k=3$.  Then 
the minimum norm of the rank 96 lattice $L(M,N,3)$ is 6 or 8.  
\end{prop}
\pf 
The value of $\mu$ depends on whether there exists rank 32 even unimodular lattices $U, M, N$ as in \refpp{gen3.5} so that $\mu (N)=4$.  
\eop

\begin{thm}\labtt{gen5} 
There exists an even unimodular lattice $L(M,N,k)$ of rank $\ell$ 
and minimum norm $\mu$ for the following 
pairs $(\ell, \mu)$: 

(i) 
$(96,8)$ (the extremal bound is 10); 
  
(ii) $(120,8)$ (the extremal bound is 12).
 
(iii) $(128,8)$ (the extremal bound is 12)
\end{thm} 
\pf 
We use \refpp{minlmn}.  

(i) Take $k=4$ and $U, M, N\cong \L$ \refpp{leechleechpolar}.  

(ii) Take $k=5$ and $U, M, N\cong \L$ \refpp{leechleechpolar}.

(iii) Take $k=4$ where $U, M, N$ are rank 32 lattices as in \refpp{gen3.5}.  
\eop

\section{Appendix: the index-determinant formula}

\begin{thm}\labtt{indexdet} {\it (``Index-determinant formula'')} 
Let $L$ be a rational lattice, and $M$ a sublattice of $L$ of finite
index $| L:M|$. Then
\begin{displaymath} det(L)\ |
L:M|^2=det(M).\end{displaymath}
\end{thm}
\pf   This is a well-known result.  
Choose a basis $x_1,\cdots, x_n$ for $L$ and positive integers $d_1,d_2
\cdots, d_n$, so that $M$ has a basis $d_1x_1,d_2x_2,\cdots,
d_nx_n$. 
A Gram matrix for the  lattice $M$ is
$G_M=((d_ix_i,d_jx_j))=DG_LD$, where
\begin{displaymath}
D = \left( \begin{array}{ccc}
 d_1 & \: & \:\\
 \: & d_2 & \:  \\
 \: & \: & \ddots
\end{array} \right),
\end{displaymath}
and $G_L=((x_i,x_j))$ is a Gram matrix for $L$.  Thus $det(G_M)=det(D)^2{\cdot} det(G_L)$.\eop

\section{Appendix: about fourvolution type sublattices and polarizations of Leech}

\begin{de}\labttr{fourvolution}  A fourvolution $f$ is a linear transformation whose square is $-1$.  If $f$ is orthogonal, $f-1$ doubles norms.  
\end{de}  

\begin{de} \labttr{fourvolutiontype}  
Let $L$ be an integral lattice.  
A sublattice $M$ of $L$ is of {\it fourvolution type} if there exists a fourvolution $f$ so that 
$M=L (f-1)$ (whence $M\cong \sqrt 2 L$).  The same terminology applies to scaled copies of $\L$.  
\end{de}

\begin{lem} \labttr{leechleechpolar} If $U\cong \L$, there are polarizations of $\U$ by sublattices $M\cong N\cong \L$.   
\end{lem}
\pf 
Here is one proof.  
We use a fact about $O(\L )$, that there are pairs of fourvolutions $f, g$ so that 
$\la f, g \ra$ is a double cover of a dihedral group of order $2k$ for which an element of odd order $k>1$ has no eigenvalue 1 on $\L$.  
There exist examples of this for $k=3, 5$, at least (for which $C_{O(\L)}(\la f, g \ra)\cong 2{\cdot}G_2(4), 2{\cdot}HJ$, respectively) \cite{gr12}.  
We take $M:=\L (f-1)$ and $N:=\L (g-1)$.  
Since $2\L=\L (f-1)^2=\L (g-1)^2$, 
$M \cap N \ge 2\L$.  
We argue that the pair $M, N$ gives a polarization.  Since $(M\cap N)/2\L$ consists of vectors fixed by $\la f, g \ra$, it is 0.  By determinant considerations, $M+N=\U$.  
\eop


\begin{thebibliography}{GRC99}


\bibitem{bachocnebe} Christine Bachoc and Gabriele Nebe, 
    Extremal lattices of minimum 8 related to the Mathieu group $M_{22}$, 
    J. Reine Angew. Math. 494 (1998) 
    

\bibitem{borcherds} Richard Borcherds, The Leech lattice, Proc. Royal Soc. London A398 (1985) 365-376. 


    
\bibitem{conway} John Conway, 
A characterization of Leech's lattice, Invent. Math. 7 (1969), 137-142.  
    
\bibitem{cs} John Conway and Neil Sloane, Sphere Packings, Lattices and Groups, Springer-Verlag 1988.  

    


\bibitem{dmn} 
%MR1650629 (99k:17048)
Dong, C.; Li, H.; Mason, G.; Norton, S. P., 
Associative subalgebras of the Griess algebra and related topics. The Monster and Lie algebras (Columbus, OH, 1996), 27--42,
Ohio State Univ. Math. Res. Inst. Publ., 7, de Gruyter, Berlin, 1998. 
\bibitem{gr12} Robert L. Griess, Jr.,   Twelve Sporadic Groups,
Springer Verlag, 1998.  




\bibitem{poe} Robert L. Griess, Jr., Pieces of Eight, Advances in Mathematics, 148, 75-104 (1999).



\bibitem{uniqe8} Robert L. Griess, Jr.,  Positive definite lattices of rank at
most 8,  Journal of Number Theory, 103 (2003),
77-84.  



\bibitem{ibw1}  Robert L. Griess, Jr., 
Involutions on the the Barnes-Wall lattices and their fixed point sublattices, I.  
Pure and Applied Mathematics Quarterly, vol.1, no. 4, (Special Issue: In Memory of Armand Borel, Part 3 of 3) 989-1022, 2005.  


\bibitem{bwy} Robert L. Griess, Jr., 
Pieces of $2^d$: existence and uniqueness for
Barnes-Wall and Ypsilanti lattices.    Advances in Mathematics, 196 (2005) 147-192.   math.GR/0403480

\bibitem{bwycorr}  Robert L. Griess, Jr., 
 Corrections and additions to `` Pieces of $2^d$: existence and uniqueness for
Barnes-Wall and Ypsilanti lattices. '' [Adv. Math. 196 (2005) 147-192], Advances in Mathematics 211 (2007) 819-824.  

\bibitem{leech} John Leech, 
Notes on sphere packings, Canadian Journal of Mathematics 19 (1967), 251-267.  

\bibitem{lm} James Lepowsky and Arne Meurman, An $E_8$ approach to the Leech lattice and the Conway groups, J. Algebra 77 (1982), 484-504.  

\bibitem{sloanesite}  Gabriele Nebe and Neil J. A. Sloane,  Catalogue of Lattices, 
http://www.research.att.com/${\sim }$njas/lattices/abbrev.html 


\bibitem{niem} H. V. Niemeier, 
Definite Quadratische Formen der Diskriminante 1 und Dimension 24, Doctoral Dissertation, G\"ottingen, 1968.  




\bibitem{serre} Jean-Pierre Serre, A Course in Arithmetic, Springer
Verlag, Graduate Texts in Mathematics 7, 1973.  


\bibitem{tits} Tits, Jacques:   Four Presentations of Leech's lattice, in Finite
Simple Groups, II, Proceedings of a London Math. Soc. Research
Symposium, Durham, 1978, ed. M. J. Colllins, pp, 306-307,
Adademic Press, London, New York, 1980.  




\bibitem{turyn} E. F. Assmus, Jr., H. F. Mattson, Jr. and R. J. Turyn,  Research to Develop the Algebraic Theory of Codes, Report AFCRL-67-0365, Air Force Cambridge Res. Labs., Bedford, Mass, June 1967.  




\end{thebibliography}
\end{document}